\documentclass[12pt]{article}
\usepackage{amsmath,amssymb}
\usepackage{url}
\begin{document}
\newtheorem{theorem}{Theorem}
\newtheorem{lemma}[theorem]{Lemma}
\newtheorem{corollary}[theorem]{Corollary}
\newtheorem{definition}[theorem]{Definition}
\newtheorem{example}[theorem]{Example}
\pagenumbering{roman}
\renewcommand{\thetheorem}{\thesection.\arabic{theorem}}
\renewcommand{\thelemma}{\thesection.\arabic{lemma}}
\newenvironment{proof}{\noindent{\bf{Proof.\/}}}{\hfill$\blacksquare$\vskip0.1in}
\renewcommand{\thetable}{\thesection.\arabic{table}}
\renewcommand{\thedefinition}{\thesection.\arabic{definition}}
\renewcommand{\theexample}{\thesection.\arabic{example}}
\renewcommand{\theequation}{\thesection.\arabic{equation}}
\newcommand{\mysection}[1]{\section{#1}\setcounter{equation}{0}
\setcounter{theorem}{0} \setcounter{lemma}{0}
\setcounter{definition}{0}}
\newcommand{\mrm}{\mathrm}
\newcommand{\beq}{\begin{equation}}
\newcommand{\eeq}{\end{equation}}

\newcommand{\ben}{\begin{enumerate}}
\newcommand{\een}{\end{enumerate}}

\newcommand{\hc}{\hat{c}}
\newcommand{\hd}{\hat{d}}
\newcommand{\hp}{\hat{p}}
\newcommand{\hq}{\hat{q}}
\newcommand{\hu}{\hat{u}}
\newcommand{\hv}{\hat{v}}
\newcommand{\hP}{\widehat{P}}
\newcommand{\hQ}{\widehat{Q}}
\newcommand{\hF}{\widehat{F}}
\newcommand{\hG}{\widehat{G}}
\newcommand{\hH}{\widehat{H}}
\newcommand{\tT}{\widetilde{T}}
\newcommand{\btB}{\mathbf{\widetilde{B}}^{(s)}}
\newcommand{\bA}{\mathbf{A}}
\newcommand{\bB}{\mathbf{B}}
\newcommand{\cT}{\check{T}}
\newcommand{\C}{\mbox{$\mathbb{C}$}}
\title
{\bf A Further Property of Functions in \\ Class ${\bf B}^{\boldsymbol(m)}$}

\author
{Avram Sidi\\
Computer Science Department\\
Technion - Israel Institute of Technology\\ Haifa 32000, Israel\\
E-mail:\quad  \url{asidi@cs.technion.ac.il}\\
URL:\quad    \url{http://www.cs.technion.ac.il/~asidi}}
\date{September 2015}
\bigskip\bigskip
\maketitle \thispagestyle{empty}
\newpage
\begin{abstract} We say that a function $\alpha(x)$ belongs to the set ${\bf A}^{(\gamma)}$ if it has an asymptotic expansion of the form
$\alpha(x)\sim \sum^\infty_{i=0}\alpha_ix^{\gamma-i}$ as $x\to\infty$, which can be differentiated term by term infinitely many times.
A function  $f(x)$ is in the class ${\bf B}^{(m)}$ if it satisfies a linear homogeneous differential equation of the form $f(x)=\sum^m_{k=1}p_k(x)f^{(k)}(x)$, with $p_k\in {\bf A}^{(i_k)}$, $i_k$ being integers satisfying $i_k\leq k$.
These functions have been shown to have many interesting properties,  and their integrals $\int^\infty_0 f(x)\,dx$, whether convergent or divergent,  can be evaluated very efficiently via the Levin--Sidi $D^{(m)}$-transformation. (In case of divergence, they are defined in some summability sense, such as Abel summability or Hadamard finite part or a mixture of these two.) In this note, we show that if
$f(x)$ is in ${\bf B}^{(m)}$, then so is $(f\circ g)(x)=f(g(x))$, where $g(x)>0$ for all large $x$ and $g\in {\bf A}^{(s)}$,  $s$ being a positive integer. This enlarges the scope of the $D^{(m)}$-transformation considerably to include functions of complicated arguments.
We demonstrate  the validity of our result with an application of the $D^{(3)}$ transformation to two integrals $I[f]$ and $I[f\circ g]$, for some $f\in{\bf B}^{(3)}$ and $g\in{\bf A}^{(2)}$.
\end{abstract}

\vspace{1cm} \noindent {\bf Mathematics Subject Classification
2010:} 34E05, 41A60, 65B05, 65D30.

\vspace{1cm} \noindent {\bf Keywords and expressions:}
Class ${\bf B}^{(m)}$ functions, infinite-range integrals, $D^{(m)}$ transformation,
acceleration of convergence,  Abel sum, Hadamard finite part, asymptotic expansions.

\thispagestyle{empty}
\newpage
\pagenumbering{arabic}
\section{Introduction and main result} \label{se1}
In this work, we continue our study of the properties of functions that belong to the class ${\bf B}^{(m)}$, which was introduced by Levin and Sidi in \cite{Levin:1981:TNC},
and treated further in Sidi \cite[Chapter 5]{Sidi:2003:PEM}. Specifically, we study the nature of the function $\phi(x)=f(g(x))$ given that  $f\in{\bf B}^{(m)}$. We  show that $\phi\in{\bf B}^{(m)}$ too, provided  $g(x)$ is a function that, roughly speaking, behaves like a polynomial as $x\to\infty$.   Thus, we conclude  that the class ${\bf B}^{(m)}$ is closed under variable transformations that behave polynomially at infinity.

Before recalling the definition of the class ${\bf B}^{(m)}$, we recall the definition of another function class  that was introduced and denoted ${\bf A}^{(\gamma)}$ also in  \cite{Levin:1981:TNC}. The classes ${\bf A}^{(\gamma)}$ feature prominently in the definition of the class ${\bf B}^{(m)}$, as will be clear soon.

\begin{definition}\label{defA}
A function $\alpha(x)$ belongs to the set $\mathbf{A}^{(\gamma)}$, where $\gamma$
is complex in general,  if it is
infinitely differentiable for all large $x>0$ and has a Poincar\'{e}-type
asymptotic expansion of the form
\begin{equation} \label{eq:ea1}
\alpha(x) \sim  \sum^{\infty}_{i=0} \alpha_ix^{\gamma-i}\quad \text{as $x \to \infty$},
\end{equation}
and its derivatives have Poincar\'{e}-type asymptotic expansions
obtained by differentiating that in \eqref{eq:ea1} formally term by
term.
\begin{itemize}
\item  If $\alpha_0 \neq 0$ in \eqref{eq:ea1}, then
$\alpha(x)$ is said to belong to $\mathbf{A}^{(\gamma)}$ strictly. In this case,
$\alpha(x)$ satisfies the asymptotic equality $\alpha(x) \sim   \alpha_0x^{\gamma}$ as $x \to \infty$.
\item
If the asymptotic expansion in \eqref{eq:ea1} is empty, that is,  $\alpha_i=0$ for all $i$ in \eqref{eq:ea1}, then either (i)\,$\alpha(x)\equiv0$ or (ii)\,$\alpha(x)=O(x^{-\mu})$ as $x\to\infty$ for every $\mu>0$, that is, $\alpha(x)$ tends to zero  faster than any negative power of $x$. (An example is $\alpha(x)=\exp(-cx^r)$, with $c>0$ and $r>0$.)
\item
If $\alpha(x)$ has an empty (nonempty) asymptotic expansion, we denote that by writing $\alpha\sim0$ ($\alpha\not\sim0$).
 \end{itemize}
\end{definition}

\noindent{\bf Remarks A.} The following are simple consequences of Definition \ref{defA}. We shall make use of  them  later. (For more, see
\cite[Chapter 5]{Sidi:2003:PEM}.)

\begin{enumerate}
\item[A1.]\label{a1} [{A1}]
$\mathbf{A}^{(\gamma)} \supset \mathbf{A}^{(\gamma-1)} \supset \mathbf{A} ^{(\gamma-2)}
\supset \cdots,$ so that if $\alpha \in \mathbf{A}^{(\gamma)}$, then, for any
positive integer $k,\ \alpha \in  \mathbf{A}^{(\gamma+k)}$ but not strictly.
Conversely, if $\alpha \in \mathbf{A}^{(\delta)}$ but not strictly, then
$\alpha \in \mathbf{A}^{(\delta-k)}$ strictly for a unique positive integer $k$.
\item[A2.] \label{a2}If $\alpha \in \mathbf{A}^{(\gamma)}$ strictly, then $\alpha \not \in\mathbf{A}^{(\gamma-1)}$.
\item[A3.] \label{a3}If $\alpha,\beta \in \mathbf{A}^{(\gamma)}$, then
$\alpha \pm \beta \in \mathbf{A}^{(\gamma)}$ as well.
If $\alpha \in \mathbf{A}^{(\gamma)}$ and
$\beta \in \mathbf{A}^{(\gamma+k)}$ strictly for some positive integer $k$,
then $\alpha \pm \beta \in \mathbf{A}^{(\gamma+k)}$ strictly.
\item[A4.] \label{a4}If $\alpha \in \mathbf{A}^{(\gamma)}$ and $\beta \in \mathbf{A}^{(\delta)}$, then
$\alpha \beta \in \mathbf{A}^{(\gamma+\delta)}$.
\item[A5.] \label{a5} If $\alpha \in \mathbf{A}^{(\gamma)}$ and $\beta \in
\mathbf{A}^{(\delta)}$ strictly, then $\alpha/\beta \in \mathbf{A}^{(\gamma-\delta)}$. If $\alpha \in \mathbf{A}^{(\gamma)}$ strictly and $\beta \in \mathbf{A}^{(\delta)}$ strictly, then $\alpha/\beta \in \mathbf{A}^{(\gamma-\delta)}$ strictly. (Note that these are not true if $\beta \in \mathbf{A}^{(\delta)}$, but not strictly.)
\item[A6.] \label{a6}If $\alpha \in \mathbf{A}^{(\gamma)}$ strictly,  such that $\alpha(x)>0$
for all large $x$, and we define
$\theta(x)=[\alpha(x)]^{\xi}$, then
$\theta\in \mathbf{A}^{(\gamma\xi)}$ strictly.
\item[A7.] \label{a7}
If $\alpha \in \mathbf{A}^{(\gamma)}$ strictly and $\beta\in
\mathbf{A}^{(k)}$ strictly for some positive integer $k$, such that
$\beta(x)>0$ for all large $x>0$,
and we define $\theta(x)=\alpha(\beta(x))$, then
$\theta\in \mathbf{A}^{(k\gamma)}$ strictly.
\item[A8.] \label{a8}
If $\alpha \in \mathbf{A}^{(\gamma)}$ (strictly) and $\gamma \neq 0$, then
$\alpha ' \in  \mathbf{A}^{(\gamma-1)}$ (strictly). If $\alpha \in \mathbf{A}^{(0)}$,
then $\alpha ' \in  \mathbf{A}^{(-2)}$.

\end{enumerate}

Now, by the way $\mathbf{A}^{(\gamma)}$
is defined, there may be any number of functions $\alpha(x)$ in $\mathbf{A}^{(\gamma)}$
having the same asymptotic expansion. (Concerning the uniqueness of $\alpha(x)$,  see the last paragraph
of \cite[Appendix A]{Sidi:2003:PEM}.) To avoid this, in certain places, it is
more convenient to work with  subsets $\mathbf{X}^{(\gamma)}$ of
$\mathbf{A}^{(\gamma)}$, which are defined next.

 \begin{definition}\label{defX} The subsets $\mathbf{X}^{(\gamma)}$ of
$\mathbf{A}^{(\gamma)}$ are defined for all $\gamma$ collectively
as follows:
\begin{enumerate}
\item [(i)]
A function $\alpha$ belongs to $\mathbf{X}^{(\gamma)}$ if
either $\alpha\equiv 0$ or $\alpha\in\mathbf{A}^{(\gamma-k)}$ strictly
for some nonnegative integer $k$.   Thus, $\alpha\sim0 \Leftrightarrow \alpha\equiv 0$ now.
\item [(ii)]
$\mathbf{X}^{(\gamma)}$ is closed under addition and multiplication by
scalars.
\item [(iii)]
If $\alpha\in\mathbf{X}^{(\gamma)}$ and
$\beta\in\mathbf{X}^{(\delta)}$, then $\alpha\beta\in\mathbf{X}^{(\gamma+\delta)}$;
if, in addition,  $\beta\in\mathbf{A}^{(\delta)}$ strictly, then
$\alpha/\beta\in\mathbf{X}^{(\gamma-\delta)}$.
\item [(iv)]
If $\alpha\in\mathbf{X}^{(\gamma)}$, then $\alpha'\in\mathbf{X}^{(\gamma-1)}$.
\end{enumerate}
\end{definition}

It is obvious that
no two functions in $\mathbf{X}^{(\gamma)}$ have the same
asymptotic expansion, since
if $\alpha,\beta\in\mathbf{X}^{(\gamma)}$, then either
$\alpha\equiv\beta$ or $\alpha-\beta\in\mathbf{A}^{(\gamma-k)}$ strictly
for some nonnegative integer $k$. Thus, $\mathbf{X}^{(\gamma)}$ does not
contain functions $\alpha(x)\not\equiv 0$ that satisfy
$\alpha(x)=O(x^{-\mu})$ as $x\to\infty$ for every $\mu>0$,
such as $\exp(-cx^s)$ with $c,s>0$.

Functions $\alpha(x)$ that are given as sums of series
$\sum^{\infty}_{i=0}\alpha_ix^{\gamma-i}$ that converge for all
large $x$ form a subset of $\mathbf{X}^{(\gamma)}$; obviously, such functions are of the form
$\alpha(x)=x^{\gamma}R(x)$ with $R(x)$ analytic at infinity.
Thus, $R(x)$ can be rational functions that are bounded at infinity,
for example.

We now turn to the definition of the class ${\bf B}^{(m)}$.

\begin{definition} \label{defB}
 A function $f(x)$ that is infinitely differentiable for all
large $x$ belongs to the set $\mathbf{B}^{(m)}$ if it satisfies a linear
homogeneous ordinary differential equation of order $m$ of the form
\begin{equation} \label{eq:ea2}
f(x) = \sum^m_{k=1} p_k(x) f^{(k)}(x),
\end{equation} where either $p_k\sim0$  or
$p_k\in\mathbf{A}^{(i_k)}$ strictly for some integer  $i_k\leq k$,  $1\leq k\leq m-1,$
and $p_m\in\mathbf{A}^{(i_m)}$ strictly for some integer  $i_m\leq m$.
\end{definition}

\noindent{\bf Remarks B.} The following are  consequences of Definition \ref{defA}.
They can be found in \cite{Levin:1981:TNC} and \cite[Chapter 5]{Sidi:2003:PEM}.
\begin{enumerate}
\item[B1.] \label{p1}
If $f \in \mathbf{B}^{(m)}$, then $f \in \mathbf{B}^{(\widehat{m})}$ for every $\widehat{m}>m$.
\item[B2.] \label{p2}
Consequently,
${\bf B}^{(1)}\subset{\bf B}^{(2)}\subset{\bf B}^{(3)}\subset\cdots\ .$
\item[B3.] \label{p3}
If $f \in \mathbf{B}^{(m)}$ with smallest $m$,  then the differential equation
\eqref{eq:ea2}   satisfied by $f(x)$ is unique, provided the $p_k$ are restricted (to ${\bf X}^{(k)}$ instead of $\mathbf{A}^{(k)}$) such that either $p_k\equiv0$ or $p_k\in{\bf X}^{(i_k)}$ strictly for some integer $i_k\leq k$,  $1\leq k\leq m$,
and $p_m\in{\bf X}^{(i_m)}$ strictly for some integer $i_m\leq m$.
(See \cite[p. 99, Proposition 5.1.5]{Sidi:2003:PEM}.)
\item[B4.] \label{p4} If $f\in{\bf A}^{(\gamma)}$ with $\gamma\neq 0$, then  $f \in \mathbf{B}^{(1)}$.
\item[B5.] \label{p5}
If $g_i\in \mathbf{B}^{(r_i)}$, $i=1,\ldots,\mu,$ then the following are true:\footnote{The assertions made in Remarks B5 and B6  are proved in \cite[pp. 107--109, Heuristics 5.4.1--5.4.3]{Sidi:2003:PEM} by relaxing the definition of ${\bf B}^{(m)}$ by assuming that the $p_k(x)$ are in some ${\bf A}^{(i_k)}$, with no restrictions on the integers $i_k$, and by making some additional assumptions. Examples suggest that  $f\in\mathbf{B}^{(m)}$ precisely as in Definition \ref{defB}, however.}

\begin{itemize}
\item $f=\prod^\mu_{i=1}g_i\in \mathbf{B}^{(m)}$,\quad $m\leq \prod^\mu_{i=1}r_i$.
 \item $f=\sum^\mu_{i=1}g_i\in \mathbf{B}^{(m)}$,\quad $m\leq \sum^\mu_{i=1}r_i$.
 \end{itemize}

\item[B6.]\label{p6} If $g_i\in \mathbf{B}^{(r)}$, $i=1,\ldots,\mu,$ and satisfy the same ordinary differential equation, then the following are true:
 \begin{itemize}
 \item  $f=\prod^\mu_{i=1}g_i\in \mathbf{B}^{(m)}$,\quad $m\leq \binom{r+\mu-1}{\mu}$. \\ In particular, if $g\in\mathbf{B}^{(r)}$, then
     $f=(g)^\mu\in \mathbf{B}^{(m)}$,\quad $m\leq \binom{r+\mu-1}{\mu}$.
 \item $f=\sum^\mu_{i=1}g_i\in \mathbf{B}^{(m)}$,\quad $m\leq r$.
\end{itemize}

\item[B7.]\label{p7} If $f \in \mathbf{B}^{(m)}$ and is integrable at infinity, then,
under some additional minor conditions at $x=\infty$,
\beq \label{eqj1} \int^\infty_x f(t)\,dt\sim \sum^{m-1}_{k=0}x^{\rho_k}f^{(k)}(x)\sum^\infty_{i=0}
\beta_{ki}x^{-i}\quad \text{as $x\to\infty$},\eeq
where $\rho_k$ are integers depending only on the $p_k(x)$ and satisfy
\beq \label{eqj2} \rho_k \leq \bar{\rho}_k=\bigg[\max_{\substack{k+1\leq n\leq m\\ p_n\not \sim0}}(i_n-n)\bigg]+k+1\leq k+1,\quad k=0,1,\ldots,m-1.\eeq

This result forms the basis of the $D^{(m)}$ transformation of Levin and Sidi \cite{Levin:1981:TNC}, which has proved to be an extremely efficient convergence accelerator for the computation of the integrals $\int^\infty_0f(x)\,dx$.
 \end{enumerate}

By Remarks  B1,  B2, B5, and B6, it is  clear that the classes $B^{(m)}$ contain an ever increasing number of functions, and this implies that the $D^{(m)}$ transformation is a comprehensive convergence acceleration method with ever increasing scope.

Finally, we would like to mention that most special functions that appear in scientific and engineering applications belong to one of the classes ${\bf B}^{(m)}$.

In this note, we continue our exploration of the properties of the classes
 ${\bf B}^{(m)}$. Analogously to what happens to the sum $f+g$ and the product $fg$ of two functions $f$ and $g$,  discussed in Remarks B5 and B6 above,  we wish to explore what happens to their composition.  Specifically, we address the following question: If $f(x)$ is in ${\bf B}^{(m)}$, then what can be said about $(f\circ g)(x)=f(g(x))$? Under what conditions on $g(x)$ is $f\circ g\in {\bf B}^{(\widetilde{m})}$ for some $\widetilde{m}$? We give an answer to this question in the next theorem, which is our main result. We provide a detailed proof of this theorem in the next section, where we make repeated use of Remarks A1--A8 without mentioning them. Finally, to keep the proof simpler, we replace the sets ${\bf A}^{(\gamma)}$ by their subsets ${\bf X}^{(\gamma)}$, even though the assertions of the theorem are true with the sets ${\bf A}^{(\gamma)}$.

\begin{theorem}\label{thfg}\begin{enumerate}
\item Let  $f(x)$ be a solution to the linear homogeneous differential equation of order $m$
\beq\label{eqr1}f(x) = \sum^m_{k=1} p_k(x) f^{(k)}(x), \eeq
such that   either $p_k\equiv0$ or
$p_k\in\mathbf{X}^{(i_k)}$ strictly for some integer  $i_k$, $1\leq k\leq m-1,$
and $p_m\in\mathbf{X}^{(i_m)}$ strictly for some integer  $i_m$. Let also $g\in{\bf X}^{(s)}$ strictly for some positive integer $s$, such that $\lim_{x\to\infty}g(x)=+\infty$.   Then $\phi(x)\equiv f(g(x))$ satisfies a linear homogeneous differential equation of order $m$ of the form
\beq\label{eqr2}\phi(x)=\sum^m_{k=1}\pi_k(x)\phi^{(k)}(x),\eeq  where the $\pi_k$ are determined by the $p_k$ and are such that either $\pi_k\equiv0$ or $\pi_k\in\mathbf{X}^{(r_k)}$ strictly for some integer  $r_k$, $1\leq k\leq m-1$, and
$\pi_m\in\mathbf{X}^{(r_m)}$ strictly for some integer  $r_m$.
Actually, we have
\begin{align}  &r_m=s(i_m-m)+m, \label{eqp33}\\ \intertext{and}
&r_k\leq \max\{s(i_k-k),r_{k+1}-(k+1),r_{k+2}-(k+2),\ldots,r_m-m\}+k, \notag\\ \label{eqp32}
& \text{if $\pi_k\not\equiv0$},\quad  k=m-1,m-2,\ldots,2,1.\end{align}(Note: On the right-hand side of the inequality in \eqref{eqp32},
$s(i_k-k)$ is absent when $p_k\equiv0$, and $r_n-n$ is absent when $\pi_n\equiv0$ for $n\in \{k+1,k+2,\ldots,m-1\}$.)
In addition,
\beq \label{eqp35}r_k\leq \max_{\substack{k\leq n\leq m \\ p_n\not\equiv0}}[s(i_n-n)]+k,\quad k=1,\ldots,m.\eeq
(The explicit expression for $\pi_m$ is given in \eqref{eqpim}. The rest of the $\pi_k$ are given by the recursion relation in \eqref{eqp12}.)

\item
Let $f(x)$ be in ${\bf B}^{(m)}$ and let $g(x)$ be in ${\bf X}^{(s)}$ strictly for some positive integer $s$, such that $\lim_{x\to\infty}g(x)=+\infty$.   Then $\phi(x)\equiv f(g(x))$ is also in ${\bf B}^{(m)}$.
\end{enumerate}
\end{theorem}

Clearly, this theorem expands  considerably the scope of the class ${\bf B}^{(m)}$,  hence the scope of the $D^{(m)}$ transformation, to include functions of complicated arguments, in the following sense: If  the $D^{(m)}$ transformation accelerates the convergence of the integral $\int^\infty_0 f(x)\,dx$,  it accelerates the convergence of the integral $\int^\infty_0 f(g(x))\,dx$ as well, when $g\in{\bf A}^{(s)}$ strictly for some positive integer $s$.

In Section \ref{se3},  we demonstrate  the validity of our result with an application of the $D^{(3)}$ transformation to two integrals $I[f]$ and $I[f\circ g]$, for some $f\in{\bf B}^{(3)}$ and $g\in{\bf A}^{(2)}$.

In Section \ref{se4}, we show via an example that $f\in{\bf B}^{(m)}$ with minimal $m$ does not necessarily mean that $m$ is minimal also for $\phi(x)$ even though $\phi\in{\bf B}^{(m)}$ too by Theorem \ref{thfg}.

\section{Proof of main result}\label{se2}
\setcounter{equation}{0}

\subsection{Preliminaries}
First, we note  that, being in ${\bf X}^{(s)}$, $g(x)$ has an asymptotic expansion of the form
 \beq\label{g1} g(x)\sim\sum^\infty_{n=0} g_nx^{s-n}\quad \text{as $x\to\infty$}\quad \text{and} \quad g_0>0,\eeq
from which, we also have that
\beq\label{g2}g(x)=\sum^s_{n=0} g_nx^{s-n}+O(x^{-1})\quad \text{as $x\to\infty$},\eeq meaning that $g(x)$ is a polynomial of degree $s$ or behaves like one as $x\to\infty$.

Thus, using  the notation
$$[a]_0=1\quad \text{and}\quad  [a]_i=\prod^{i-1}_{j=0}(a-j), \quad  i=1,2,\ldots,$$ we also have that
$$g^{(i)}(x)\sim\sum^\infty_{n=0} g_i[s-n]_ix^{s-n-i}\quad \text{as $x\to\infty$},\quad i=1,2,\ldots, $$
from which,
$$ g^{(i)}(x)\sim [s]_ig_0 x^{s-i}\quad \text{as $x\to\infty$},\quad i=0,1,2,\ldots,s,$$  and, since $(\sum^s_{n=0}g_nx^{s-n})^{(i)}=0$ for $i\geq s+1$,
$$g^{(i)}(x)\sim [-\mu]_ig_{s+\mu}x^{-\mu-i}\quad \text{as $x\to\infty$},\quad i=s+1,s+2,\ldots, $$ where $g_{s+\mu}$, $\mu\geq1$, is the first nonzero $g_{s+j}$ with $j\geq1$, assuming that $g^{(i)}\not\equiv0$. Of course,
$$ g_{s+j}=0, \quad j=1,2,\ldots,\quad\Rightarrow\quad g^{(i)}(x)\equiv0,\quad
i=s+1,s+2,\ldots, $$ and this can occur when $g(x)=\sum^s_{n=0}g_nx^{s-n}$, for example, in which case, $g^{(i)}(x)\equiv0$ for $i=s+1,s+2,\ldots\ .$

Summarizing, we have
\beq\label{eqp3}g^{(i)}(x)>0\  \text{for all large $x$},\quad  g^{(i)}\in {\bf X}^{(\tau_i)}\  \text{strictly},\quad \tau_i=s-i,\quad i=0,1,\ldots,s,\eeq and
\begin{align}\label{eqp3e} &\text{either}\ g^{(i)}\in {\bf X}^{(\tau_i)}\  \text{strictly},\quad \tau_i=-\mu-i<s-i, \quad i=s+1,s+2,\ldots,\notag\\ &\text{or}\quad g^{(i)}(x)\equiv0,\quad i=s+1,s+2,\ldots\ .\end{align}
Clearly, \beq\label{eqp75}\tau_i\leq s-i, \quad s=0,1,\ldots\ .\eeq

 Next, it is clear that we  need to prove that $\phi(x)=f(g(x))$ satisfies \eqref{eqr2}, with $\pi_k\in
{\bf X}^{(r_k)}$ for some integer $r_k$ when $\pi_k\not\equiv0$. Replacing $x$ by $g(x)$ throughout the differential equation \eqref{eqr1} satisfied by $f(x)$, we have
\beq \label{eqp2}f(g(x))=\sum^m_{k=1} p_k(g(x)) f^{(k)}(g(x)),\eeq and this is the starting point of our proof.
Here, we emphasize that $f^{(k)}(g(x))$ stands for the $k$th derivative of $f$ with respect to its {\em argument}, evaluated at $g(x)$; that is, $f^{(k)}(g(x))=[\frac{d^k}{dt^k}f(t)]|_{t=g(x)}$.
Thus, $f^{(k)}(g(x))$  does {\em not}  stand for the
$k$th derivative of $\phi(x)=f(g(x))$ with respect to $x$.

Whenever convenient, in  the sequel, we   write $p_k$, $\pi_k$, and $g^{(i)}$ instead of $p_k(x)$, $\pi_k(x)$, and $g^{(i)}(x)$, respectively,  for short. Thus, $p_k(g)$ and  $f^{(k)}(g)$ stand for $p_k(g(x))$ and $f^{(k)}(g(x))$, respectively.

\subsection{Special cases}
Before embarking on the proof for   arbitrary $m$, we look at the  simple but instructive cases involving $m=1,2$.
\subsubsection{The case $m=1$}
Here we consider two different cases.
\begin{itemize}\item {\em The case $f\in{\bf A}^{(\gamma)}$ strictly, $\gamma\neq0$}:\\
In this case $\phi(x)\equiv (f\circ g)(x)= f(g(x))$ satisfies the identity
$$ \phi(x)=\pi_1(x)\phi'(x),\quad \pi_1(x)=\frac{\phi(x)}{\phi'(x)}=\frac{f(g(x))}{f'(g(x))g'(x)}.$$
Because $f\in{\bf A}^{(\gamma)}$ and $\gamma\neq0$, we have that $f'\in{\bf A}^{(\gamma-1)}$ strictly. Consequently,
$f(g(x))\in{\bf A}^{(s\gamma)}$ and  $f'(g(x))\in{\bf A}^{(s(\gamma-1))}$.
This implies that $\pi_1\in {\bf A}^{(r_1)}$ strictly, where
$$ r_1=s\gamma-[s(\gamma-1)+(s-1)]=1.$$
Thus, $(f\circ g) \in{\bf B}^{(1)}$.
\item {\em The general case of $f(x)=p_1(x)f'(x)$}:\\
In this case, $f(x)=p_1(x)f'(x)$,  $p_1\in {\bf A}^{(i_1)}$ strictly, with $i_1$ an integer.
Replacing $x$ by $g(x)$, this differential equation becomes
$$ f(g(x))=p_1(g(x))f'(g(x))\quad \Rightarrow \quad p_1(g(x))=\frac{f(g(x))}{f'(g(x))}.$$
Now, $\phi(x)\equiv (f\circ g)(x)= f(g(x))$ satisfies
$$ \phi(x)=\pi_1(x)\phi'(x),\quad \pi_1(x)=\frac{\phi(x)}{\phi'(x)}=\frac{f(g(x))}{f'(g(x))g'(x)}=\frac{p_1(g(x))}{g'(x)}.$$
Therefore, by the fact that $p_1(g(x))\in {\bf A}^{(si_1)}$ strictly, we have that $\pi_1\in {\bf A}^{(r_1)}$ strictly, where
$$ r_1=si_1-(s-1)=s(i_1-1)+1.$$ Now, when $i_1\leq 1$, we have that $f\in{\bf B}^{(1)}$. In this case,
$r_1\leq 1$,   which implies that $(f\circ g) \in{\bf B}^{(1)}$ too.
\end{itemize}
\bigskip
\subsubsection{The case $m=2$}
 With $\phi(x)=f(g(x))$, we have
 $$ \phi'(x)=f'(g(x))g'(x),\quad \phi''(x)=f''(g(x))(g'(x))^2+f'(g(x))g''(x).$$
Substituting these in \eqref{eqr2},  we obtain
$$ f(g(x))=\pi_1(x)[f'(g(x))g'(x)]+\pi_2(x)[f''(g(x))(g'(x))^2+f'(g(x))g''(x)],$$ which, upon rearranging, becomes
$$f(g(x))=[\pi_1(x)g'(x)+\pi_2(x)g''(x)]f'(g(x))+[\pi_2(x)(g'(x))^2]f''(g(x)).$$
Comparing this with \eqref{eqp2}, we identify
the following equations for $\pi_1$ and $\pi_2$:
\begin{align*} p_1(g(x))&=\pi_1(x)g'(x)+\pi_2(x)g''(x)\\ p_2(g(x))&=\pi_2(x)(g'(x))^2. \end{align*}

Since $s>0$,  $g'\in{\bf X}^{(s-1)}$  strictly, and positive for all large $x$. Therefore, $$\pi_2=\frac{p_2(g)}{(g')^2}\in{\bf X}^{(r_2)}\ \text{strictly}, \quad r_2=si_2-2(s-1)=s(i_2-2)+2.$$ Next,
\begin{align*} \pi_1&=\frac{p_1(g)}{g'}-\frac{\pi_2g''}{g'}\in{\bf X}^{(r_1)}\ \text{strictly if $\pi_1\not\equiv0$},\\
\intertext{and since ${p_1(g)}/{g'}\in{\bf X}^{(si_1-(s-1))}$ strictly and
$\pi_2g''/g'\in{\bf X}^{(r_2+(s-2)-(s-1))}$, we also have}
 r_1&\leq\max\{si_1-(s-1), r_2+(s-2)-(s-1)\} \\ &=
\max\{s(i_1-1), r_2-2\}+1 \\ &=\max\{s(i_1-1),s(i_2-2)\}+1.
\end{align*}
Note that if $p_1\equiv0$, the term $s(i_1-1)$ is absent throughout.

\subsection{The case of arbitrary $m$}
We prove part 1 of the theorem first. We start with \eqref{eqp2}.
By the Fa\`{a} di Bruno formula for differentiation of $f\circ g$, we have

\beq\label{eqp4} \phi^{(n)}(x)=\frac{d^n}{dx^n}f(g(x))=\sum^n_{k=1}
B_{n,k}(g'(x),g''(x),\ldots,g^{(n-k+1)}(x)) f^{(k)}(g(x)),\eeq where
 $ B_{n,k}(y_1,y_2,\ldots,y_{n-k+1})$ is the Bell polynomial defined as in
\beq\label{eqp5a} B_{n,k}(y_1,y_2,\ldots,y_{n-k+1})=\sum\frac{n!}{\prod^{n-k+1}_{i=1}(j_i!)}
\prod^{n-k+1}_{i=1}\bigg(\frac{y_i}{i!}\bigg)^{j_i},\eeq the summation being on the nonnegative integers $j_1,j_2,\ldots, j_{n-k+1}$ such that
\beq\label{eqp5c}\sum^{n-k+1}_{i=1}j_i=k\quad \text{and}\quad \sum^{n-k+1}_{i=1}i j_i=n.\eeq
(The simplest of these polynomials are
$ B_{n,1}(y_1,\ldots,y_n)=y_n$ {and} $B_{n,n}(y_1)=y_1^n.$)

For the Fa\`{a} di Bruno formula, see Johnson \cite{Johnson:2002:CHF}, for example.
For Bell polynomials, see Bell \cite{Bell:1934:EP}.

Now, if the conjecture in \eqref{eqr2} is true, then substituting \eqref{eqp4} in \eqref{eqr2}, we must have
\beq\label{eqp99} f(g(x))=\sum^m_{n=1}\pi_n(x)\bigg[\sum^n_{k=1}
B_{n,k}(g'(x),g''(x),\ldots,g^{(n-k+1)}(x)) f^{(k)}(g(x))\bigg],\eeq
which, upon rearrangement, becomes
\beq\label{eqp5} f(g(x))=\sum^m_{k=1}\bigg[\sum^m_{n=k}\pi_n(x)
B_{n,k}(g'(x),g''(x),\ldots,g^{(n-k+1)}(x))\bigg] f^{(k)}(g(x)).\eeq
Comparing \eqref{eqp5} with \eqref{eqp2}, we realize that the equalities
\beq\label{eqp7} p_k(g(x))=\sum^m_{n=k}\pi_n(x)
B_{n,k}(g'(x),g''(x),\ldots,g^{(n-k+1)}(x)),\quad k=1,\ldots,m,\eeq must hold.
Clearly, this is an $m$-dimensional upper  triangular system of linear equations for $\pi_1(x),\ldots,\pi_m(x)$, provided the latter exist.
The diagonal of the matrix of this system is $[B_{1,1}(g'(x)),B_{2,2}(g'(x)),\ldots,
B_{m,m}(g'(x))]$, and, by the fact that $B_{n,n}(y_1)=y_1^n$,  $$[B_{1,1}(g'(x)),B_{2,2}(g'(x)),\ldots,
B_{m,m}(g'(x))]=[(g'(x))^1,(g'(x))^2,\ldots,(g'(x))^m].$$
Since $g'(x)>0$ for all large $x$, this diagonal is positive, hence the linear system in \eqref{eqp7} has  a unique solution for $\pi_1(x),\ldots,\pi_m(x)$. With the existence of the $\pi_k$ established, we now need to show that, $\pi_k\in{\bf X}^{(r_k)}$  strictly for some integer $r_k$ when $\pi_k\not\equiv0$. We achieve this goal by induction on $k$, in the order $k=m,m-1,\ldots,2,1.$

To be able to proceed, we need to analyze  $B_{n,k}(g',g'',\ldots,g^{(n-k+1)})$.
By \eqref{eqp5a},
\beq \label{eqp57} B_{n,k}(g',g'',\ldots,g^{(n-k+1)})=\sum\frac{n!}{\prod^{n-k+1}_{i=1}(j_i!)}
\prod^{n-k+1}_{i=1}\bigg(\frac{g^{(i)}}{i!}\bigg)^{j_i},\eeq
the summation being on the nonnegative integers $j_1,j_2,\ldots, j_{n-k+1}$ subject to the constraints in \eqref{eqp5c},
 and by \eqref{eqp3} and
\eqref{eqp3e}, when  $g^{(i)}\not\equiv0$, $1\leq i\leq n-k+1$,
$$ \prod^{n-k+1}_{i=1}(g^{(i)})^{j_i}\in{\bf X}^{(\sigma(j_1,\ldots,j_{n-k+1}))},\quad \sigma(j_1,\ldots,j_{n-k+1})=
\sum^{n-k+1}_{i=1}\tau_ij_i.$$ Upon invoking  \eqref{eqp75}, this gives
$$\sigma(j_1,\ldots,j_{n-k+1})\leq \overline{\sigma}(j_1,\ldots,j_{n-k+1})=\sum^{n-k+1}_{i=1}(s-i)j_i.$$
Of course,  this also means that
$$ \prod^{n-k+1}_{i=1}(g^{(i)})^{j_i}\in{\bf X}^{(\overline{\sigma}(j_1,\ldots,j_{n-k+1}))}\begin{cases} \text{strictly},\quad \text{if\ $n-k+1\leq s$},\quad \text{by \eqref{eqp3}}\\ \text{not strictly},\quad\text{otherwise, \quad by \eqref{eqp3e}.}\end{cases}$$
At first sight, $ \overline{\sigma}(j_1,\ldots,j_{n-k+1})$ seems to depend on $j_1,\ldots,j_{n-k+1}$. This is not so, however.
In fact,  on account of the constraints in \eqref{eqp5c}, $ \overline{\sigma}(j_1,\ldots,j_{n-k+1})$ depends only on $n$ and $k$:
$$ \overline{\sigma}(j_1,\ldots,j_{n-k+1})=s\sum^{n-k+1}_{i=1}j_i-\sum^{n-k+1}_{i=1}i j_i=sk-n.$$
Consequently, because all the terms in the summation on the right-hand side of \eqref{eqp57} are  in ${\bf X}^{(sk-n)}$, we have
\beq\label{eqp11} L_{n,k}\equiv B_{n,k}(g',g'',\ldots,g^{(n-k+1)})\in {\bf X}^{(sk-n)},\quad\text{but not necessarily strictly.} \footnote{Note that,
in our study, we need only be concerned with $g^{(i)}$ for $i=1,\ldots,m$.   (i)\,When  $s<m$, because of \eqref{eqp3e}, not all products $\prod^{n-k+1}_{i=1}(g^{(i)})^{j_i}$ are in ${\bf X}^{(\overline{\sigma}(j_1,\ldots,j_{n-k+1}))}$ strictly (some of them may even be zero identically), and this results in \eqref{eqp11}. (ii)\,When  $s\geq m$, however,
on account of \eqref{eqp3}, the products  $\prod^{n-k+1}_{i=1}(g^{(i)})^{j_i}$ are all positive and in ${\bf X}^{(\overline{\sigma}(j_1,\ldots,j_{n-k+1}))}$ strictly. Therefore, $L_{n,k}$ are all positive and in ${\bf X}^{(sk-n)}$ strictly.}\eeq
 By the fact that $B_{k,k}(g')=(g')^k$, however, we have
\beq\label{eqv1} L_{k,k}\in {\bf X}^{(sk-k)}\ \text{strictly}.\eeq

We now start the induction with $\pi_m$, which we obtain from the last of the equations in \eqref{eqp7}. Thus,
\beq\label{eqpim}\pi_m B_{m,m}(g')=p_m(g)\quad \Rightarrow \quad \pi_m=\frac{p_m(g)}{B_{m,m}(g')}=\frac{p_m(g)}{(g')^m}. \eeq
By the fact that $p_m(g)\in{\bf X}^{(si_m)}$ strictly and $(g')^m\in{\bf X}^{((s-1)m)}$ strictly,
it is clear that $\pi_m\not\equiv0$ and
\beq\label{eqp22} \pi_m\in {\bf X}^{(r_m)}\ \text{strictly},\quad r_m=si_m-m(s-1)=s(i_m-m)+m.\eeq
We have thus shown the validity of our assertion for $\pi_m$.

We now continue by induction on $k$.  Let us assume  that the assertion is true also for $\pi_{m-1},\pi_{m-2},\ldots,\pi_{k+1}$, namely,  $\pi_n\in{\bf X}^{(r_n)}$ strictly for some integer $r_n$ if $\pi_n\not\equiv0$,  $n= m-1,m-2,\ldots,k+2,k+1.$  The proof will be complete if we show that $\pi_k\in{\bf X}^{(r_k)}$ strictly for some integer $r_k$ if $\pi_k\not\equiv0$.  Solving \eqref{eqp7}, namely, $p_k(g)=\sum^m_{n=k}\pi_nL_{n,k}$,  for $\pi_k$, we obtain
\beq \label{eqp12}\pi_k=\frac{p_k(g)}{L_{k,k}}-\sum^m_{n=k+1}\pi_n\frac{L_{n,k}}{L_{k,k}}.\eeq

First, $L_{k,k}\in {\bf X}^{(sk-k)}$ {strictly}, and  $p_k(g)\in{\bf X}^{(si_k)}$ strictly when $p_k\not\equiv0$; therefore, if $p_k\not\equiv0$,
\beq\label{eqp14} \frac{p_k(g)}{L_{k,k}}\in{\bf X}^{(\mu_k)} \ \text{strictly},\quad \mu_k=si_k-(sk-k)=s(i_k-k)+k.\eeq If $p_k\equiv0$, then $p_k(g)\equiv0$ too, and, therefore,
${p_k(g)}/{L_{k,k}}\equiv0.$

Next, for $n=k+1,k+2,\ldots,m,$ by \eqref{eqp11} and \eqref{eqv1} and the induction hypothesis, if $\pi_n\not\equiv0$,
\beq\label{eqp15}  \pi_n\frac{L_{n,k}}{L_{k,k}}\in {\bf X}^{(\nu_{n,k})} \ \text{strictly},\quad \nu_{n,k}\leq  r_n+(sk-n)-(sk-k)=(r_n-n)+k.\eeq
If $\pi_n\equiv0$, then $\pi_nL_{n,k}/L_{k,k}\equiv0$ too.

Combining \eqref{eqp14} and \eqref{eqp15} in  \eqref{eqp12},  when $\pi_k\not\equiv0$, we have
\beq\label{eqp16}  \pi_k\in{\bf X}^{(r_k)}\ \text{strictly},\quad r_k\leq \max\{\mu_k,\nu_{k+1,k},\nu_{k+2,k},\ldots,\nu_{m,k}\}.\footnote{On the right-hand side of the equality for $r_k$ in \eqref{eqp16},
$\mu_k$ is absent when $p_k\equiv0$, and $r_n-n$ is absent when $\pi_n\equiv0$ for $n\in \{k+1,k+2,\ldots,m-1\}$. }\eeq
This completes the induction step.

By the fact that $\mu_k-k=s(i_k-k)$ and $\nu_{n,k}-k\leq r_n-n$, the   equality for $r_k$ in \eqref{eqp16} is   identical to that in  \eqref{eqp32}.
Finally, the proof of the inequality in \eqref{eqp35} can now be achieved by induction on $k$, in the order $k=m,m-1,\ldots,2,1.$
This completes the proof of part 1 of the theorem.

Part 2 is a straightforward corollary of part 1, since $f\in {\bf B}^{(m)}$ means that $f(x)$ is exactly as in part 1, with $i_k\leq k$ (equivalently,  $i_k-k\leq0$),  for each $k$. Invoking this in \eqref{eqp33}--\eqref{eqp35}, the proof of part 2 is completed.

\section{The $D^{(m)}$ transformation and an application}\label{se3}
\setcounter{equation}{0}
\subsection{The $D^{(m)}$ transformation}
The $D^{(m)}$ transformation for computing infinite-range  integrals of the form
 $I[f]=\int^\infty_0f(t)\,dt$ is defined as follows:
\begin{itemize}
\item Choose a sequence $\{x_l\}^\infty_{l=-1}$, such that
$$ 0=x_{-1}<x_0<x_1<x_2<\cdots,\quad \lim_{l\to\infty}x_l=\infty.$$
\item Define $F(x)=\int^x_0f(t)\,dt$ and compute $F(x_l)$, $l=0,1,\ldots\ .$
This is best achieved by computing the integrals $\chi_i=\int^{x_i}_{x_{i-1}}f(t)\,dt,$
$i=0,1,\ldots,$ numerically  (preferably by a low order Gaussian quadrature formula),  and forming  $F(x_l)=\sum^l_{i=0}\chi_i.$
\item Let $n=(n_1,n_2,\ldots,n_{m})$, where $n_1,n_2,\ldots,n_{m}$ are nonnegative integers, and solve the $(N+1)$-dimensional linear system
\beq \label{eqh1} F(x_l)=D^{(m,j)}_n+\sum^{m}_{k=1}x_l^{\rho_{k-1}} f^{(k-1)}(x_l)\sum^{n_k-1}_{i=0}\frac{\bar{\beta}_{ki}}{x_l^i},\quad j\leq l\leq j+N, \quad N=\sum^m_{k=1}n_k,\eeq
for $D^{(m,j)}_n$, which is the approximation to $I[f]$. Here $\rho_k$ are as in \eqref{eqj1} and \eqref{eqj2}, and
$\bar{\beta}_{k,i}$ are additional auxiliary unknowns, which are not of interest.

In case the $\rho_k$ are not known but  the $i_k$ are known, the $\rho_k$ in
 \eqref{eqh1} can be replaced  by their upper bounds $\bar{\rho}_k$ given in \eqref{eqj2}. If the $i_k$ too are not known, we can replace the $\rho_k$ by their  ultimate upper bounds $k+1$ given again in \eqref{eqj2}. Thus, the user-friendly version of the $D^{(m)}$ transformation is now defined as in
\beq  \label{eqh2}F(x_l)=D^{(m,j)}_n+\sum^{m}_{k=1}x_l^{k} f^{(k-1)}(x_l)\sum^{n_k-1}_{i=0}\frac{\bar{\beta}_{ki}}{x_l^i},\quad j\leq l\leq j+N, \quad N=\sum^m_{k=1}n_k.\eeq

\end{itemize}

From the way the user-friendly version of the $D^{(m)}$ transformation is defined
as in \eqref{eqh2}, it is clear that we need not know anything about the differential equation satisfied by $f(x)$.   First, an upper bound for the order of the differential equation can be taken as $m$, as suggested by Remarks B5 and B6. (In most cases, we can determine the smallest $m$ quite easily.) Next, we need to be able to compute $f^{(i)}(x)$, $i=0,1,\ldots,m-1$. Finally,
since the $x_l$ are at our disposal, we  can be  choose them   appropriately to ensure excellent  convergence rates.

When $f\in {\bf B}^{(m)}$, the sequences of approximations $\{ D^{(m,j)}_{(\nu,\nu,\ldots,\nu)}\}^\infty_{\nu=0}$, with fixed $j$ (in particular, with $j=0$), have the best convergence properties. These sequences can be computed very efficiently by applying
the W$^{(m)}$ algorithm of Ford and Sidi \cite{Ford:1987:AGR}.\footnote{The W$^{(m)}$
algorithm, when used for implementing the user-friendly   $D^{(m)}$ transformation defined via the linear systems in \eqref{eqh2},
 is designed to compute the sequences $\{A^{(j)}_N\}^\infty_{N=0}$ recursively, via the solutions of $$ F(x_l)= A^{(j)}_N+\sum^m_{k=1}x_l^kf^{(k-1)}(x_l)\sum^{\lfloor{(N-k)/m}\rfloor}_{i=0}\frac{\bar{\beta}_{ki}}{x_l^i},\quad j\leq l\leq j+N.$$ Then $\{D^{(m,j)}_{(\nu,\nu,\ldots,\nu)}\}^\infty_{\nu=0}$ is a proper subsequence of $\{A^{(j)}_N\}^\infty_{N=0}$. In fact, $A^{(j)}_{m\nu}=D^{(m,j)}_{(\nu,\nu,\ldots,\nu)},$ $\nu=0,1,\ldots,$ with $A^{(j)}_0=D^{(m,j)}_{(0,0,\ldots,0)}=F(x_j).$
See \cite[Section 7.3, p. 165]{Sidi:2003:PEM}.}
(For $m=1$, the W$^{(m)}$ algorithm reduces to the W algorithm of Sidi \cite{Sidi:1982:ASC}.)

Note that, for determining $D^{(m,j)}_n$, we need as input the finite-range integrals $F(x_l)$, $j\leq l\leq j+N.$ In case $I[f]$ exists as a regular improper integral,
$F(x_{j+N})$ is naturally the best available approximation to $I[f]$ out of the integrals
$F(x_l)$, $j\leq l\leq j+N$. Therefore,  it is instructive to compare the accuracy of
$D^{(m,j)}_n$ with that of $F(x_{j+N})$. It is always observed that the accuracy of $D^{(m,j)}_n$ is much higher  than  that of $F(x_{j+N})$. (See the application that follows.)

\subsection{An application}
We would like to apply the $D^{(m)}$ transformation to the  integrals $I[f]$ and
$I[\phi]$ (with $g(x)=x^2$), where
$$ f(x)=\frac{\sin^2 x}{x^2}\quad \text{and}\quad
\phi(x)=f(x^2)=\frac{\sin^2 x^2}{x^4}.$$
We have (see Gradshteyn and Ryzhik \cite[formulas 3.821.9 and 3.852.3]{Gradshteyn:1980:TIS})
$$ I[f]=\frac{\pi}{2}\quad \text{and}\quad I[\phi]=\frac{2\sqrt{\pi}}{3}.$$

Now $f\in{\bf B}^{(3)}$ since
$$f(x)=\sum^3_{k=1}p_k(x)f^{(k)}(x);\quad p_1(x)=-\frac{2x^2+3}{4x},\ \ p_2(x)=-\frac{3}{4},\ \ p_3(x)=-\frac{x}{8},$$ with $p_k\in{\bf A}^{(i_k)}$ strictly,
$$ i_1=1,\quad i_2=0,\quad i_3=1.$$
(Note that $p_k\in{\bf X}^{(i_k)}$ strictly as well.)

By Theorem \ref{thfg},
$\phi\in{\bf B}^{(3)}$ too, because
$$\phi(x)=\sum^3_{k=1}\pi_k(x)f^{(k)}(x);\quad \pi_1(x)=-\frac{16x^4+15}{64x^3},\ \ \pi_2(x)=-\frac{9}{64x^2},\ \ \pi_3(x)=-\frac{1}{64x},$$
with $\pi_k\in{\bf A}^{(r_k)}$ strictly,
$$ r_1=1,\quad r_2=-2,\quad r_3=-1.$$ Letting $s=2$ in Theorem \ref{thfg}, we see that the $r_k$ are consistent with \eqref{eqp33}--\eqref{eqp35}.
The $\pi_k$ are obtained from
$$\pi_3(g')^3={p_3(g)},\quad
\pi_2{(g')^2}+3\pi_3g'g''={p_2(g)},\quad
 \pi_1g'+\pi_2g''+\pi_3g'''=p_1(g).$$
(Note that $\pi_k\in{\bf X}^{(r_k)}$ strictly as well.)

 We have applied the user-friendly version of the $D^{(m)}$ transformation (as defined in \eqref{eqh2})  to $I[f]$ with $x_l=1.6(l+1)$ and to $I[\phi]$ with $x_l=\sqrt{1.6(l+1)}$, $l=0,1,\ldots\ .$ The results of these computations are given in Table \ref{table1}. In this table,  we compare the errors in $D^{(3,0)}_{(\nu,\nu,\nu)}[f]$ with the corresponding errors in $F(x_{3\nu})$. Similarly, we compare
the errors in $D^{(3,0)}_{(\nu,\nu,\nu)}[\phi]$ with the corresponding errors in $\Phi(x_{3\nu})$. (As before, $F(x)=\int^x_0f(t)\,dt$ and, similarly,  $\Phi(x)=\int^x_0\phi(t)\,dt$.)
This comparison demonstrates very clearly the power  of the $D^{(m)}$ transformation as a convergence accelerator.

 \begin{table}[htbp]
\caption{\label{table1}  Results from the $D^{(3)}$ transformation applied to (i)\,$I[f]=\int^\infty_0f(x)\,dx$ with $x_l=1.6(l+1)$ and (ii)\,$I[\phi]=\int^\infty_0\phi(x)\,dx$ with $x_l=\sqrt{1.6(l+1)}$, where $f(x)=(\sin x/x)^2$ and $\phi(x)=f(x^2)$. We have defined $F(x)=\int^x_0f(t)\,dt$ and $\Phi(x)=\int^x_0\phi(t)\,dt$.}\vspace{-1cm}
\begin{center}
$$
\begin{array}{||r||c|c||c|c||}
\hline
\nu &|F(x_{3\nu})-I[f]|&|D^{(3,0)}_{(\nu,\nu,\nu)}[f]-I[f]|&|\Phi(x_{3\nu})-I[\phi]|&|D^{(3,0)}_{(\nu,\nu,\nu)}[\phi]-I[\phi]|\\
\hline\hline
0&       3.44D-01&       3.44D-01&   9.64D-02&       9.64D-02  \\
1&       7.86D-02&       7.06D-02&   1.03D-02&       7.06D-03  \\
2&       4.40D-02&       6.96D-03&   4.36D-03&       8.89D-04  \\
3&       3.17D-02&       1.69D-04&   2.66D-03&       8.63D-07  \\
4&       2.37D-02&       5.70D-07&  1.72D-03&       1.88D-06  \\
5&       1.98D-02&       2.48D-07&  1.32D-03&       5.73D-08  \\
6&       1.62D-02&       1.32D-08&  9.73D-04&       9.57D-09  \\
7&       1.44D-02&       1.04D-10&  8.14D-04&       4.62D-11  \\
8&       1.23D-02&       7.15D-11&  6.47D-04&       1.65D-10  \\
9&       1.13D-02&       3.88D-12&  5.65D-04&       1.09D-11  \\
10&       9.98D-03&       4.83D-13&   4.70D-04&       3.58D-13 \\
\hline
     \end{array}
$$
\end{center}
 \end{table}
\section{A further development}\label{se4}
While reviewing the class ${\bf B}^{(m)}$ in Section \ref{se1}, we mentioned that  if the degree $m$ of the differential equation \eqref{eq:ea2} is minimal, then the differential equation is unique. In view of Theorem \ref{thfg}, we might be led to think that
$m$ is also the minimal degree of the differential equation \eqref{eqr2} satisfied by $\phi(x)=f(g(x))$. In other words, it might sound plausible that if $m$ is the smallest integer for which $f\in{\bf B}^{(m)}$, then $m$ is the  smallest integer for which   $\phi\in{\bf B}^{(m)}$ as well. We show via an example that this is not always the case; that is, it is possible that $\phi\in{\bf B}^{(\widehat{m})}$ for some $\widehat{m}< m$. Note that Theorem \ref{thfg} does not contradict this since, by Remark B1,   $\phi\in{\bf B}^{(\widehat{m})}$ implies $\phi\in{\bf B}^{(\mu)}$ for every  $\mu>\widehat{m}$, thus for $\mu=m$ in particular.

\sloppypar
Consider the function $f(x)=1/(\sqrt{x}+1)^{3}$. Now
$f\not\in{\bf B}^{(1)}$. If $f(x)$ were in ${\bf B}^{(1)}$, then we would have  $f(x)=p_1(x)f'(x)$ with $p_1\in {\bf A}^{(i_1)}$ for some integer $i_1\leq 1$. But
$$ p_1(x)=\frac{f(x)}{f'(x)}=-\frac{2}{3}(x+\sqrt{x})\not\in{\bf A}^{(\gamma)}
\quad \text{for any $\gamma$.}$$
It is true, however, that $f\in{\bf B}^{(2)}$. To see this, we observe that
$$ f(x)=\bigg[\frac{ \sqrt{x}-1}{x-1}\bigg]^3= f_1(x)+f_2(x),$$ where $$f_1(x)=\frac{\sqrt{x}(x+3)}{(x-1)^3},\quad f_2(x)=-\frac{3x+1}{(x-1)^3}. $$
 Now, $f_1\in{\bf A}^{(-3/2)}$, while $f_2\in {\bf A}^{(-2)}$, hence, by Remark B4,
$f_1\in{\bf B}^{(1)}$ and $f_2\in{\bf B}^{(1)}$. By Remark B5, $f=f_1+f_2\in {\bf B}^{(2)}$ since we have already seen that $f\not\in{\bf B}^{(1)}$.

Let us now turn to $\phi(x)=f(g(x))$ with $g(x)=x^2$. We have
$\phi(x)=1/(x+1)^3$. Clearly, $\phi\in{\bf A}^{(-3)}$, hence $\phi\in{\bf B}^{(1)}$ by Remark B4.


\begin{thebibliography}{1}

\bibitem{Bell:1934:EP}
E.T. Bell.
\newblock Exponential polynomials.
\newblock {\em Ann. Math.}, 35:258--277, 1934.

\bibitem{Ford:1987:AGR}
W.F. Ford and A.~Sidi.
\newblock An algorithm for a generalization of the {Richardson} extrapolation
  process.
\newblock {\em SIAM J. Numer. Anal.}, 24:1212--1232, 1987.

\bibitem{Gradshteyn:1980:TIS}
I.S. Gradshteyn and I.M. Ryzhik.
\newblock {\em {Table of Integrals, Series, and Products}}.
\newblock Academic Press, New York, 1980.
\newblock Forth printing 1983.

\bibitem{Johnson:2002:CHF}
W.P. Johnson.
\newblock The curious history of {Fa\`{a} di Bruno's} formula.
\newblock {\em Amer. Math. Monthly}, 109:217--234, 2002.

\bibitem{Levin:1981:TNC}
D.~Levin and A.~Sidi.
\newblock Two new classes of nonlinear transformations for accelerating the
  convergence of infinite integrals and series.
\newblock {\em Appl. Math. Comp.}, 9:175--215, 1981.
\newblock Originally appeared as a Tel Aviv University preprint in 1975.

\bibitem{Sidi:1982:ASC}
A.~Sidi.
\newblock An algorithm for a special case of a generalization of the
  {Richardson} extrapolation process.
\newblock {\em Numer. Math.}, 38:299--307, 1982.

\bibitem{Sidi:2003:PEM}
A.~Sidi.
\newblock {\em {Practical Extrapolation Methods: Theory and Applications}}.
\newblock Number~10 in Cambridge Monographs on Applied and Computational
  Mathematics. Cambridge University Press, Cambridge, 2003.

\end{thebibliography}

\end{document}